\input amstex
\UseAMSsymbols
\input epsf.tex
\documentstyle{amsppt}


\loadbold \nologo \pageheight{8.5truein} \pagewidth{7.0truein}
\topmatter
\title A counterexample to the second inequality of Corollary (19.10)in
the monograph "Ricci Flow and the Poincare Conjecture" by J.Morgan
and G.Tian
\endtitle
\author Abbas Bahri
\endauthor
\endtopmatter

\medskip
\subhead {1.Introduction}\endsubhead
\medskip

We provide in the sequel a counterexample to the second inequality
of Corollary (19.10) of [1], in line with what was announced in [2].
\medskip
\subhead {2. The counterexample} \endsubhead
\medskip

Observe that the inequality implies that a curve which is a geodesic
for a given metric remains a geodesic for the metric evolved through
the Ricci flow as the curve itself is evolved through the curve
shortening flow $H$.

However, the norm of $H$ is $k$ and therefore the curve itself does
not move under the curve shortening flow.

It follows that, if the inequality holds, a curve that is a geodesic
remains, without moving, a geodesic for the evolved metric through
the Ricci flow. We provide below a counterexample to this
conclusion.

The equation of a geodesic reads:

$$\ddot x^k+\Gamma_{ij}^k \dot x^i \dot x^j=0$$

We introduce coordinates along a small piece of curve of this
geodesic so that this small piece of curve defines the $x^1$-axis of
coordinates. It follows that, for this set of coordinates, on this
small piece of curve:

$$\dot x^1(t) \gneq 0$$, whereas $\dot x^s(t)=0, s\neq 1$

Then, on this small piece of curve, the geodesic equation verified
for $s \neq 1$ reads:

$$\ddot x^s+\Gamma_{11}^s \dot x^1 \dot x^1=0   \Longleftrightarrow \Gamma_{11}^s \dot x^1 \dot x^1=0$$

As the metric evolves through the Ricci flow, $\ddot x^s$ remains
equal to zero, for $s \neq 1$, on the piece of curve.
$\Gamma_{11}^s$ was zero at the time zero of the evolution. It
follows that the first variation of $\Gamma_{11}^s$, which we denote
$\delta \Gamma_{11}^s$ should also be zero along this piece of
curve.

We choose an arbitrary point $x_0$ on this piece of curve. We may
assume that, at this point, the metric tensor reads
$g_{ij}(x_0)=\delta_{ij}$. There is no loss of generality in this
requirement.

Then,
$$\Gamma_{11}^s(x_0)=1/2(2 \frac {\partial g_{1s}}{\partial x^1}-\frac {\partial g_{11}}{\partial
x^s})(x_0)$$

The formula for $\delta \Gamma_{11}^s$ reads (use geodesic normal
coordinates in order to derive the formula. $\delta \Gamma_{ij}^k$
are components of a tensor. Therefore, the computation in any set of
coordinates, eg geodesic normal coordinates, provides the formula
for another set of coordinates):

$$\delta\Gamma_{11}^s(x_0)=1/2(2 \nabla_1\delta g_{1s}-\nabla_s \delta g_{11})(x_0)$$

where we used the fact that $g_{ij}(x_0)=\delta_{ij}$.

We know that $\delta g_{ij}=-2R_{ij}$, so that:

$$\delta\Gamma_{11}^s(x_0)=(\nabla_s R_{11}-2\nabla_1
R_{1s})(x_0)=0$$

We have:

$$R_{11}=\frac {\partial \Gamma_{11}^t}{\partial x^t}-\frac {\partial \Gamma_{1t}^t}{\partial
x^1}+O(\Gamma^2)$$

$$R_{1s}=\frac {\partial \Gamma_{1s}^t}{\partial x^t}-\frac {\partial \Gamma_{1t}^t}{\partial
x^s}+O(\Gamma^2)$$

The expression of the covariant derivative of $R_{11}$ and $R_{1s}$
is complicated; but there are third derivatives of the metric tensor
in  $\delta \Gamma_{11}^s (x_0)$. They are derived as if we were
computing in a geodesic normal coordinates system. They are the
third derivatives in:

$$\frac {\partial}{\partial x^s}(\frac {\partial \Gamma_{11}^t}{\partial x^t}-\frac {\partial \Gamma_{1t}^t}{\partial
x^1})(x_0)-2\frac {\partial}{\partial x^1}(\frac {\partial
\Gamma_{1s}^t}{\partial x^t}-\frac {\partial \Gamma_{1t}^t}{\partial
x^s})(x_0)$$

Computing, we find:
$$-(2\frac {\partial ^3 g_{st}}{\partial x^{1 ^2} \partial x^t}-2\frac {\partial ^3 g_{1s}}{\partial x^{t ^2} \partial
x^1}+\frac {\partial ^3 g_{11}}{\partial x^{t ^2} \partial
x^1}-\frac {\partial ^3 g_{tt}}{\partial x^{1 ^2} \partial
x^s})(x_0)$$

For $t \neq s, t\neq 1$, we claim that $(\frac {\partial ^3
g_{tt}}{\partial x^{1 ^2} \partial x^s})(x_0)$, for $t\neq s, s\neq
1$, is a free parameter along this piece of curve:  indeed, the only
coordinate that is non-zero along this piece of curve is $x^1$.
Thus, all $\dot x^r, r\neq 1$ are zero and the geodesic equation
involves only the Christoffel symbols $\Gamma^m_{11}$. $g_{tt}$, $t
\neq 1, t\neq s$, can interfere which this equation, but only
through the coefficients $g^{mq}$ of the inverse of the metric
tensor in front of the Christoffel symbols. These $g^{mq}$s can be
subject to constraints along the piece of curve. But, no derivative
of these $g^{mq}$, hence no derivative of $g_{tt}$, $t$ as above, is
involved in the equations that are verified. All the terms with
derivatives of the metric in these Christoffel symbols involve
$g_{1q}$, ie one of the indexes is $1$.

Thus no transversal derivative for $g_{tt}$, $t$ as above, is
involved in the geodesic equation. The behavior of $g_{tt}$ and its
derivatives along $x^1$ transversally to the piece of curve is
completely free. The conclusion follows

The other terms involve only products and powers of the first and
second derivatives of the metric tensor, on this piece of curve.
Also, on this piece of curve, only $\partial/{\partial x^1}$ is a
derivative along the curve. The other derivatives
$\partial/{\partial x^t}, t \neq 1$ are transverse to this piece of
curve.

($\frac {\partial ^3 g_{tt}}{\partial x^{1 ^2} \partial x^s})(x_0),
t \neq 1, t \neq s$ can be taken non-zero and large as we please in
a small neighborhood of the curve, so that
$\delta\Gamma_{11}^s(x_0)$ is non-zero, a contradiction.
\newpage

\widestnumber\no{99999}

\font\tt cmr12 at 24 truept 


\noindent\Refs\nofrills{\tt References}


\tenpoint

\medskip
\ref\no 1\by J.Morgan and G.Tian\book Ricci Flow and the Poincare
Conjecture\publ Clay Mathematics Monograph, AMS and Clay
Institute\vol 3\yr 2007\endref
\medskip
\ref\no 2\by A.Bahri\paper Five gaps in Mathematics \jour Advanced
Non-linear Studies\vol Vol. 15, No. 2 \yr 2015\pages 289-320\endref
 \end{document}